\documentclass[11pt]{article}
\usepackage{latexsym, amsfonts}
\newcommand{\be}{\begin{equation}}
\newcommand{\ee}{\end{equation}}
\newcommand{\bea}{\begin{eqnarray}}
\newcommand{\eea}{\end{eqnarray}}
\newcommand{\bean}{\begin{eqnarray*}}
\newcommand{\eean}{\end{eqnarray*}}
\newcommand{\brray}{\begin{array}}
\newcommand{\erray}{\end{array}}
\newcommand{\ben}{\begin{equation}{nonumber}}
\newcommand{\een}{\end{equation}{nonumber}}

\newtheorem{dfn}{Definition}[section]
\newtheorem{thm}[dfn]{Theorem}
\newtheorem{lmma}[dfn]{Lemma}
\newtheorem{ppsn}[dfn]{Proposition}
\newtheorem{crlre}[dfn]{Corollary}
\newtheorem{xmpl}[dfn]{Example}
\newtheorem{rmrk}[dfn]{Remark}
\newcommand{\bdfn}{\begin{dfn}}
\newcommand{\bthm}{\begin{thm}}
\newcommand{\blmma}{\begin{lmma}}
\newcommand{\bppsn}{\begin{ppsn}}
\newcommand{\bcrlre}{\begin{crlre}}
\newcommand{\bxmpl}{\begin{xmpl}}
\newcommand{\brmrk}{\begin{rmrk}}
\newcommand{\edfn}{\end{dfn}}
\newcommand{\ethm}{\end{thm}}
\newcommand{\elmma}{\end{lmma}}
\newcommand{\eppsn}{\end{ppsn}}
\newcommand{\ecrlre}{\end{crlre}}
\newcommand{\exmpl}{\end{xmpl}}
\newcommand{\ermrk}{\end{rmrk}}

\newcommand{\IC}{\mathbb{C}}

\newcommand{\IR}{\mathbb{R}}

\newcommand{\IT}{\mathbb{T}}

\newcommand{\IZ}{\mathbb{Z}}

\newcommand{\al}{\alpha}

\newcommand{\gma}{\gamma}

\newcommand{\cla}{{\cal A}}
\newcommand{\clb}{{\cal B}}
\newcommand{\clc}{{\cal C}}

\newcommand{\clf}{{\cal F}}
\newcommand{\clg}{{\cal G}}
\newcommand{\clh}{{\cal H}}
\newcommand{\cli}{{\cal I}}

\newcommand{\cll}{{\cal L}}

\newcommand{\cls}{{\cal S}}

\newcommand{\clu}{{\cal U}}
\newcommand{\clv}{{\cal V}}

\def\a*{{\cal A}_{h,*}}
\def\B{{\cal B}(h)}
\def\B1{{\cal B}_1(h)}
\def\b{{\cal B}^{\rm s.a.}(h)}
\def\b1{{\cal B}^{\rm s.a.}_1(h)}

\newcommand{\ot}{\otimes}

\newcommand{\raro}{\rightarrow}

\newcommand{\lgl}{\langle}
\newcommand{\rgl}{\rangle}

\def \qed {$\Box$}

\begin{document}

\begin{center}
{\large {\bf Quantum  Group of Isometries in Classical and Noncommutative Geometry}}\\
by\\
{\large Debashish Goswami {\footnote {The author gratefully acknowldges support obtained  from the Indian National Academy of Sciences through the grants for a project on `Noncommutative Geometry and Quantum Groups', and also wishes to thank    The Abdus
Salam ICTP (Trieste), where a major part of the work  was done
during a visit as Junior Assciate.}}}\\
{\large Stat-Math Unit, Kolkata Centre,}\\
{\large Indian Statistical Institute}\\
{\large 203, B. T. Road, Kolkata 700 108, India}\\
\end{center}
\begin{abstract}
We formulate a quantum generalization of the notion of the group of Riemannian isometries for a compact Riemannian manifold, by introducing a
natural notion of smooth and isometric action by a compact quantum group on a classical or noncommutative manifold described by spectral triples,
and then proving the existence of a universal object (called the quantum isometry group) in the category of compact quantum groups acting smoothly
 and isometrically on a given (possibly noncommutative) manifold satisfying certain regularity assumptions. In fact, we  identify the quantum isometry group with the universal object in a bigger category, namely the category of  `quantum families of smooth isometries', defined along the line of Woronowicz and Soltan. We also construct a spectral triple on the Hilbert space of forms on a noncommutative manifold which is equivariant  with respect to a natural unitary representation of the quantum isometry group.  
 We give explicit description of quantum isometry groups of commutative and noncommutative tori, and in this context,
 obtain the quantum double torus defined in \cite{hajac} as the
 universal quantum group of holomorphic isometries of the
 noncommutative torus.
\end{abstract}
\section{Introduction}
Since the formulation of quantum automorphism groups by Wang
(\cite{free}, \cite{wang}), following suggestions of Alain Connes,
many interesting examples of such quantum groups, particularly the
quantum permutation groups of finite sets and finite graphs, have
been extensively studied by a number of mathematicians (see, e.g.
\cite{ban1}, \cite{ban2}, \cite{univ1} and references therein),
who have also found applications to and interaction with areas
like free probability  and subfactor theory. The underlying basic
 principle of defining a quantum automorphism group corresponding to some given mathematical structure (for example, a finite set, a graph,
  a $C^*$ or von Neumann algebra) consists of two steps : first, to identify (if possible) the group of automorphisms of the structure as a universal object in a
   suitable category, and then, try to look for the universal object in a similar but bigger category by replacing groups by quantum groups of appropriate type.
    However, most of the work done so far concern some kind of quantum automorphism groups of a `finite' structure, for example, of finite sets or finite dimensional
     matrix algebras. It is thus quite natural to try to extend these ideas to the
     `infinite' or `continuous' mathematical structures, for
     example classical and noncommutative manifolds. In the
     present article, we have made an attempt to formulate and
     study the quantum analogues of the groups of Riemannian
     isometries, which play a very important role in the classical
     differential geometry. The group of Riemannian isometries of a compact Riemannian manifold $M$ can be viewed as the universal object in the
      category of all compact metrizable groups acting on $M$, with smooth and isometric action. Therefore, to define the quantum isometry group,
       it is reasonable to  consider a
     category of compact quantum groups which act on the manifold
     (or more generally, on a noncommutative manifold given by
     spectral triple) in a `nice' way, preserving the Riemannian
     structure in some suitable sense, to be precisely formulated.
     In this article, we have given a definition of such `smooth
     and isometric' action by a compact quantum group on a
     (possibly noncommutative) manifold, extending the notion of  smooth and
     isometric action by a group on a classical manifold. Indeed, the meaning of isometric action is nothing but that the action should commute
      with the `Laplacian' coming from the spectral triple, and we should mention that this idea was  already present in \cite{ban2},
      though only in the context of  a finite metric space or a finite graph.
       The
     universal object in the category of such quantum groups, if
     it
     exists, should be thought of as the quantum analogue of the
     group of isometries, and we have been able to prove its
     existence under some regularity assumptions, all of which can
     be verified for a general compact connected Riemannian
     manifold as well as the standard examples of noncommutative
     manifolds. 
       Motivated by the ideas of Woronowicz and Soltan, we actually consider a bigger category. The isometry group of a classical manifold, viewed as a compact metrizable space (forgetiing the group structure), can be seen to be the universal object of  a category whose object-class consists of subsets (not necessarily subgroups) of the set of smooth isometries of the manifold. Then it can be proved that this universal compact set has a canonical group structure. A natural quantum analogue of this has been formulated by us, called the category of `quantum families of smooth isometries'.  The underlying $C^*$-algebra of the quantum isometry group has been identified with its universal object and moreover, it is shown to be equipped with a canonical coproduct making it into a compact quantum group.  
     
     We believe that a detailed study of quantum
      isometry groups will not only give many new and interesting
      examples of compact quantum groups, it will also contribute
      to the understanding of quantum group covariant spectral
      triples. In fact, we have made some progress in this direction already by constructing a spectral triple (which is often closely related to the original spectral triple) on the Hilbert space of forms which is equivriant with respect to a canonical unitary representation of the quantum isometry group.  
      
      In a companion article \cite{jyotish} with J. Bhowmick, we  provide explicit computations of quantum
     isometry groups of a few classical and noncommutative
     manifolds.  However, we
      briefly quote some of main results of \cite{jyotish} in the
      present article. One interesting observation is that the
      quantum isometry group of the noncommutative two-torus
      $\cla_\theta$ (with the canonical spectral triple) is
       (as a $C^*$ algebra) a direct sum of two commutative
      and two noncommutative tori, and contains as a quantum subgroup (which is universal for certain class of isometric actions called holomorphic isometries)
       the  `quantum double-torus' discovered and studied by Hajac and Masuda
       (\cite{hajac}).
\section{Definition of the quantum isometry group}
\subsection{Isometry groups of classical manifolds}
We begin with a well-known characterization of the isometry group
of a (classical) compact Riemannian manifold. Let $(M,g)$ be a
compact Riemannian manifold and let $\Omega^1= \Omega^1(M)$ be the
space of smooth one-forms, which has a right
Hilbert-$C^\infty(M)$-module structure given by the
$C^\infty(M)$-valued inner product $<< \cdot, \cdot >>$ defined by
$$<< \omega, \eta >>(m)=<\omega(m),\eta(m)>|_m,$$ where
$<\cdot,\cdot>|_m$ is the Riemannian metric on the cotangent space
$T^*_mM$ at the point $m \in M$. The Riemannian volume form allows
us to make $\Omega^1$ a pre-Hilbert space, and we denote its
completion by $\clh_1$. Let $\clh_0=L^2(M,{\rm dvol})$ and consider the
de-Rham differential $d$ as an unbounded linear map from $\clh_0$
to $\clh_1$, with the natural domain $C^\infty(M) \subset \clh_0$,
and also denote its closure by $d$. Let $\cll:=- d^*
d$. The following identity can be verified by direct and easy
computation using the local coordinates : $$ (\partial \cll)(\phi,
\psi) \equiv \cll(\bar{\phi}
\psi)-\cll(\bar{\phi})\psi-\bar{\phi}\cll(\psi)=
 2<< d \phi, d \psi >>~~~{\rm for ~} \phi, \psi \in C^\infty(M)~~~~(*).$$
 \begin{ppsn}
 \label{iso}
 A smooth map $\gma : M \raro M$ is a Riemannian isometry if and only if $\gma$ commutes with $\cll$ in the sense that $\cll(f \circ \gma)=(\cll(f)) \circ \gma $ for all $f \in C^\infty(M)$.
 \eppsn
 {\it Proof :}\\
 If $\gma$ commutes with $\cll$ then from the identity (*) we get for $m \in M$ and $\phi, \psi \in C^\infty(M)$ :
 \bean
 \lefteqn{<d\phi|_{\gma(m)}, d\psi|_{\gma(m)}>|_{\gma(m)}}\\
 &=& << d \phi, d \psi >>(\gma(m))\\
 &=& \frac{1}{2}(\partial \cll (\phi, \psi) \circ \gma)(m)\\
 &=& \frac{1}{2}\partial \cll (\phi \circ \gma, \psi \circ \gma)(m)\\
 &=& << d(\phi \circ \gma), d(\psi \circ \gma) >> (m)\\
  &=& < d(\phi \circ \gma)|_m, d(\psi \circ \gma)|_m>|_m\\
  &=& <(d\gma|_m)^*(d\phi|_{\gma(m)}), (d\gma|_m)^*(d \psi|_{\gma(m)})>|_m,\\
  \eean
  which proves that $(d\gma|_m)^* : T^*_{\gma(m)}M \raro T^*_m M$ is an isometry. Thus, $\gma$ is a Riemannian  isometry.

  Conversely, if $\gma$ is an isometry, both the maps induced by $\gma$ on $\clh_0$ and $\clh_1$, i.e.
  $U^0_\gma : \clh_0 \raro \clh_0$ given by $U^0_\gma(f)=f \circ \gma$ and $U^1_\gma : \clh^1 \raro \clh^1$ given by
  $U^1_\gma(f d \phi)=(f \circ \gma) d(\phi \circ \gma)$  are unitaries. Moreover, $d \circ U^0_\gma =U^1_\gma \circ  d$ on $C^\infty(M)
  \subset \clh_0$. From this, it follows that $\cll=-d^*d $  commutes with $U^0_\gma$.
  \qed

 Now let us consider a compact metrizable (i.e. second countable) space $Y$ with a continuous map $\theta : M \times Y \raro M$. We abbreviate $\theta(m,y)$ as $ym$ and denote by $\xi_y$ the map $M \ni m \mapsto ym$.  Let $\alpha: C(M) \raro C(M) \ot C(Y) \cong C(M \times Y)$ be the map given by $\alpha(f)(m,y):=f(ym)$ for $ y \in Y$, $m \in M$ and $f \in C(M)$. For a state $\phi$ on $C(Y)$, denote by $\alpha_\phi$ the map $({\rm id} \ot \phi) \circ \alpha : C(M) \raro C(M)$. We shall also denote by $\clc$ the subspace of $C(M) \ot C(Y)$ generated by elements of the form $\alpha(f)(1 \ot \psi)$, $f \in C(M), \psi \in C(Y)$. Since $C(M)$ and $C(Y)$ are commutative algebras, it is easy to see that $\clc$ is a $\ast$-subalgebra of $C(M) \ot C(Y)$.  Then we have the following
  \bthm
  \label{classic}
  (i) $\clc$ is norm-dense in $C(M) \ot C(Y)$ if and only if for every $y \in Y$, $\xi_y$ is one-to-one.\\
  (ii) The map $\xi_y$ is $C^\infty$ for every $y \in Y$ if and only if $\alpha_\phi (C^\infty(M)) \subseteq C^\infty(M)$
  for all $\phi$. \\
  (iii) Under the hypothesis of (ii),  each $\xi_y$ is also an isometry if and only if $\alpha_\phi$ commutes with $(\cll-\lambda)^{-1}$ for all state $\phi$ and all $\lambda$ in the resolvent of $\cll$ (equivalently, $\alpha_\phi$ commutes with the Laplacian $\cll$ on $C^\infty(M)$).
  \ethm
  {\it Proof :}\\
  (i) First, assume that $\xi_y$ is one-to-one for all $y$. By Stone-Weirstrass Theorem, it is enough to show that $\clc$ separates points. Take $(m_1,y_1) \ne (m_2,y_2)$ in $M \times Y$. If $y_1 \ne y_2$, we can choose $\psi \in C(Y)$ which separates $y_1$ and $y_2$, hence $(1 \ot \psi) \in \clc$ separates $(m_1,y_1)$ and $(m_2,y_2)$. So, we can consider the case  when $y_1=y_2=y$ (say), but $m_1 \ne m_2$. By injectivity of $\xi_y$, we have $ym_1 \ne ym_2$, so there exists $f \in C(M)$ such that $f(ym_1) \ne f(ym_2)$, i.e. $\alpha(f)(m_1,y) \ne \alpha(f)(m_2,y)$. This proves the density of $\clc$.
  
  For the converse, we argue as in the proof of Proposition 3.3 of \cite{vandaelenotes}.  Assume that $\clc$ is dense in $C(M) \ot C(Y)$, and let $y \in Y$, $m_1,m_2 \in M$ such that $ym_1=ym_2$. That is, $\alpha(f)(1 \ot \psi)(m_1,y)=\alpha(f)(1 \ot \psi)(m_2,y)$ for all $f \in C(M)$, $\psi \in C(Y)$. By the density of $\clc$ we get $\chi(m_1,y)=\chi(m_2,y)$ for all $\chi \in C(M \times Y)$, so $(m_1,y)=(m_2,y)$, i.e. $m_1=m_2$.
  
  (ii) The `if part' of (ii) follows by considering the states corresponding to point evaluation, i.e. $C(Y) \ni \psi \mapsto \psi(y)$, $y \in Y$. For the converse, we note that an arbitrary state $\phi$ corresponds to a regular Borel measure $\mu$ on $Y$ so that $\phi(h)=\int h d \mu$, and thus, $\alpha_\phi(f)(m)=\int f(ym) d \mu(y)$ for $f \in C(M)$. From this, by interchanging differentiation and integation (which is allowed by the Dominated Convergence Theorem, since $\mu$ is a finite measure) we can prove that $\alpha_\phi(f)$ is $C^\infty$ whenever $f$ is so. 
  
  The assertion (iii) follows from Proposition \ref{iso} in a straghtforward way.
  \qed \vspace{1mm}\\
Let us recall a few well-known facts about the Laplacian $\cll$, viewed as a negative self-adjoint operator on the Hilbert space $L^2(M, {\rm dvol})$. It is known (see \cite{rosenberg} and references therein) that $\cll$ has compact resolvents and all its eigenvectors belong to $C^\infty(M)$. Moreover, it follows from the Sobolev Embedding Theorem that $$ \bigcap_{n \geq 1} {\rm Dom}(\cll^n)=C^\infty(M).$$ Let $\{ e_{ij}, j=1,..., d_i; i=1,2,... \}$ be the set of (normalised) eigenvectors of $\cll$, where  $e_{ij} \in C^\infty(M)$ is an eigenvector corresponding to the eigenvalue $\lambda_i$, $|\lambda_1| < |\lambda_2| <...$. We have the following:
\blmma
The complex linear span of $\{ e_{ij} \}$ is norm-dense in $C(M)$.
 \elmma {\it Proof :}\\
This is  a
consequence of the asymptotic estimates of eigenvalues
$\lambda_i$, as well as the uniform bound of the eigenfunctions
$e_{ij}$. For example, it is known (\cite{lapest},Theorem 1.2)
that there exist constants $C,C^\prime$ such that $\| e_{ij}
\|_\infty \leq C | \lambda_i|^{\frac{n-1}{4}},$ $d_i \leq C^\prime
|\lambda_i |^{\frac{n-1}{2}},$ where $n$ is the dimension of the
manifold $M$. Now, for $f \in C^\infty(M) \subseteq \bigcap_{k
\geq 1}{\rm Dom}(\cll^k)$, we write $f$ as an a-priori
$L^2$-convergent series $\sum_{ij} f_{ij} e_{ij}$ ($f_{ij} \in
\IC$), and observe that $\sum |f_{ij}|^2 |\lambda_i|^{2k} <
\infty$ for every $k \geq 1$. Choose and fix sufficiently large
$k$ such that $\sum_{i \geq 0} |\lambda_i|^{n-1-2k} < \infty$,
which is possible due to  the well-known Weyl asymptotics of
eigenvalues of $\cll$. Now,  by the Cauchy-Schwarz inequality and the estimate
for $d_i$, we have  $$ \sum_{ij} |f_{ij}| \| e_{ij} \|_\infty \leq
C (C^\prime)^{\frac{1}{2}} \left( \sum_{ij} |f_{ij}|^2 |
\lambda_i|^{2k} \right)^{\frac{1}{2}} \left(\sum_{i \geq 0}
|\lambda_i|^{n-1-2k} \right)^{\frac{1}{2}}< \infty.$$ Thus, the
series $\sum_{ij} f_{ij} e_{ij}$ converges to $f$ in sup-norm, so
${\rm Sp}\{ e_{ij},j=1,2,...,d_i;i=1,2,... \}$ is dense in sup-norm in $C^\infty(M)$, hence in $C(M)$ as
well. \qed

 Let us denote ${\rm Sp}\{ e_{ij},j=1,...,d_i;i \geq 1 \}$ by $\cla^\infty_0$ from now on. We shall now show that $C^\infty(M)$ can be replaced by the smaller subspace $\cla^\infty_0$ in Theorem \ref{classic}. We need a lemma for this,  which will be useful later on too.
 \blmma
 \label{tracepres_1}
 Let $\clh_1,\clh_2$ be Hilbert spaces and for $i=1,2$, let $\cll_i$ be (possibly unbounded) self-adjoint operator on $\clh_i$ with compact resolvents, and let   $\clv_i$ be the linear span of eigenvectors of $\cll_i$. Moreover, assume that there is an eigenvalue of $\cll_i$ for which the eigenspace is one-dimensional, say spanned by a unit vector $\xi_i$. Let $\Psi$ be a linear map from $\clv_1$ to $\clv_2$ such that $\cll_2 \Psi=\Psi \cll_1$ and $\Psi(\xi_1)=\xi_2.$ Then we have \be \label{etc} \lgl \xi_2, \Psi(x)\rgl=\lgl \xi_1, x \rgl~~\forall x \in \clv_1. \ee
 \elmma
 {\it Proof:}\\
 By hypothesis on $\Psi$, it is clear that there is a common eigenvalue, say $\lambda_0$, of $\cll_1$ and $\cll_2$, with the eigenvectors $\xi_1$ and $\xi_2$ respectively. Let us write the set of eigenvalues of $\cll_i$ as a disjoint union $\{ \lambda_0 \} \bigcup \Lambda_i$ ($i=1,2$), and let the corresponding orthogonal decomposition of $\clv_i$ be given by $\clv_i=\IC \xi_i \bigoplus_{\lambda \in \Lambda_i} \clv^\lambda_i \equiv \IC \xi_i \oplus \clv^\prime_i$, say, where $\clv^\lambda_i$ denotes the eigenspace of $\cll_i$ corresponding to the eigenvalue $\lambda$. By assumption, $\Psi$ maps $\clv^\lambda_1$ to $\clv^\lambda_2$ whenever $\lambda$ is an eigenvalue of $\cll_2$, i.e. $\clv^\lambda_2 \ne \{ 0 \}$, and otherwise it maps $\clv^\lambda_1$ into $\{ 0 \}$. Thus, $\Psi(\clv_1^\prime) \subseteq \clv^\prime_2$. Now, (\ref{etc}) is obviously satisfied for $x=\xi_1$, so it is enough to prove (\ref{etc}) for all $x \in \clv^\prime_1$. But we have $\lgl \xi, x \rgl=0$  for $x \in \clv^\prime_1$, and  since $\Psi(x) \in \clv^\prime_2=\clv_2 \bigcap \{ \xi_2 \}^\perp$, it follows that $\lgl \xi_2, \Psi(x) \rgl=0=\lgl \xi_1, x \rgl$. 
 \qed
 
      \blmma
      \label{classic2}
     Let $Y$ and $\alpha$ be as in Theorem \ref{classic}. Then the following are equivalent.\\
     (a) For every $y \in Y$, $\xi_y$ is smooth isometric.\\
     (b) For every state $\phi$ on $C(Y)$, we have $\alpha_\phi(\cla^\infty_0) \subseteq \cla^\infty_0$,  and $\alpha_\phi \cll=\cll \alpha_\phi$  on $\cla^\infty_0$.\\
      \elmma
     
      {\it Proof:}\\
    We prove only the nontrivial implication $(b) \Rightarrow (a)$.  Assume  (b) that
     $\alpha_\phi$ leaves  $\cla^\infty_0$ invariant and commutes with $\cll$ on it, for 
     every state $\phi$.
      To prove that $\alpha$ is a smooth isometric action, it is enough 
     (see the proof of  Theorem \ref{classic}) to prove that  $\alpha_y(\cla^\infty) \subseteq \cla^\infty$ for all $y \in Y$, where    $\alpha_y(f):=({\rm id} \ot {\rm ev}_y)(f)=f \circ \xi_y$, ${\rm ev}_y$ being the evaluation at the point $y$. 
     Let $M_1,...,M_k$ be the connected components of the compact manifold $M$. Thus, the Hilbert space $L^2(M, {\rm dvol})$ admits an orthogonal decomposition $\oplus_{i=1}^k L^2(M_i,{\rm dvol})$,  
     and the Laplacian $\cll$ is of the form $\oplus_i \cll_i$ where $\cll_i$ denotes the Laplacian on $M_i$. Since each $M_i$ is connected, we have ${\rm Ker}(\cll_i)=\IC \chi_i$, where $\chi_i$ is the constant function on $M_i$ equal to $1$. Now, we note that for fixed $y$ and $i$, the image of $M_i$ under the  continuous function $\xi_y$ must be mapped into a  component, say $M_j$.  Thus, by applying Lemma \ref{tracepres_1} with $\clh_1=L^2(M_i)$,$\clh_2=L^2(M_j)$, $\Psi=\xi_y$ and the $L^2$-continuity of the map $f \mapsto \alpha_y(f)=f \circ \xi_y$, we have $$ \int_{M_j} \alpha_y(f)(x) {\rm dvol}(x)=\int_{M_i} f(x) {\rm dvol}(x)$$ for all $f$ in the linear span of eigenvectors of $\cll_i$, hence (by density) for all $f$ in $  L^2(M_i)$. It follows that $\int_M \alpha_y(f) {\rm dvol}=\int_M f {\rm dvol}$ for all $f \in L^2(M)$,  in particular for all $f \in C(M)$. Since $\alpha_y$ is a $\ast$-homomorphism on $C(M)$, we have $$ \lgl \alpha_y(f), \alpha_y(g) \rgl=\int_M \alpha_y(\overline{f}g) {\rm dvol}=\int_M \overline{f}g {\rm dvol}=\lgl f, g \rgl,$$ for all $f,g \in C(M)$. Thus, $\alpha_y$ extends to an isometry on $L^2(M)$, to be denoted by the same notation, which by our assumption commutes with    the self-adjoint operator $\cll$ on the core $\cla^\infty_0$, and hence  $\alpha_y$ commutes with $\cll^n$ for all $n$.  In particular it leaves invariant the domains of each $\cll^n$, which implies   $\alpha_y(\cla^\infty) \subseteq \cla^\infty$.
   \qed

    In view of the fact that the set of isometries of $M$, denoted by $ISO(M)$, is a compact second countable (i.e. compact metrizable)  group, we see that $ISO(M)$ is the maximal compact second countable   group acting on $M$ such that the action is smooth and isometric. In other words, if we consider a catogory whose objects are compact metrizable groups acting smoothly and isometrically on $M$, and morphisms are the group homomorphisms commuting with the actions  on $M$, then $ISO(M)$ (with its canonical action on $M$) is the initial object of this cateogory. However, one can take a more general viewpoint and consider the category of compact metrizable spaces $Y$ equipped with a continuous map $\theta : M \times Y \raro M$ satisfying (i)-(iii) of Theorem \ref{classic}, or equivalently, the pair of commutative unital $C^*$-algebras $\clb=C(Y)$ and a unital  $C^*$-homomorphism $\alpha: C(M) \raro C(M) \raro \clb$ satisfying the conditions (i)-(iii).  The set of isometries $ISO(M)$ (as a topological space)  can be identified with the universal object of this category, and then one can prove that it has a group structure.
    
 It is  quite natural to formulate a quantum analogue of the above, by considering, in the spirit of Woronowicz and Soltan (see \cite{woro_pseudo} and \cite{soltan}), `quantum families of isometries', which can be defined to be a pair $(\clb, \alpha)$ where $\clb$ is a (not necessarily commutative) $C^*$-algebra and $\alpha : C(M) \raro C(M) \ot \clb$ is unital $C^*$-homomorhism satisfying (i)-(iii) of Theorem \ref{classic}, i.e. the linear span of $\alpha(C(M))(1 \ot \clb)$ (which is not necessarily a $\ast$-subalgebra any more, $\clb$ being possibly noncommutative) is norm-dense in $C(M) \ot \clb$ and for every state $\phi$ on $\clb$, the map $\alpha_\phi$ keeps $C^\infty(M)$ invariant and commutes with the Laplacian $\cll$. The morphisms of this category are obvious. We shall prove that this category has a universal object, and this universal object can be equipped with a canonical quantum group structure. This will define the quantum isometry group of a manifold.  However, we shall go beyond classical manifolds and  define quantum isometry
 group $QISO({\cla^\infty},\clh,D)$ for a spectral triple $({\cla^\infty}, \clh, D)$, with ${\cla^\infty}$ being unital, and satisfying certain assumptions. To this end, we need to carefully formulate the notion of Laplacian in noncommutative  geometry, which is the goal of the next subsection.
 
 \subsection{Laplacian in noncommutative geometry}
  Given a spectral triple $({\cla^\infty},\clh,D)$,
  we recall from \cite{fro} and \cite{con}
  the construction of the space of one-forms. We have a derivation from ${\cla^\infty}$ to the ${\cla^\infty}$-${\cla^\infty}$ bimodule $\clb(\clh)$ given by $a \mapsto [D,a]$.
  This induces  a bimodule morphism $\pi$ from $\Omega^1({\cla^\infty})$ (the bimodule of universal one-forms on ${\cla^\infty}$)  to $\clb(\clh)$,  such that
  $\pi(\delta(a))=[D,a]$, where $\delta :{\cla^\infty} \raro \Omega^1({\cla^\infty})$ denotes the universal derivation map.
  We set $\Omega^1_D \equiv \Omega^1_D({\cla^\infty}):=\Omega^1({\cla^\infty})/{\rm Ker}(\pi) \cong \pi(\Omega^1({\cla^\infty}))\subseteq \clb(\clh)$.
  Assume that the spectral  triple is of compact type and has a finite dimension in the sense of Connes (\cite{con}),
  i.e. there is some $p>0$ such that the operator $|D|^{-p}$ (interpreted as the inverse of the restriction of $|D|^p$ on the closure of its range,
  which has a finite co-dimension since $D$ has compact resolvents) has finite nonzero Dixmier trace,
  denoted by $Tr_\omega$ (where $\omega$ is some suitable Banach limit,  see, e.g. \cite{con}, \cite{fro}).
  Consider the canonical `volume form' $\tau $ coming from the Dixmier trace, i.e.
  $\tau : \clb(\clh) \raro \IC$ defined by $\tau(A):=\frac{1}{Tr_\omega(|D|^{-p})} Tr_\omega(A |D|^{-p}).$ Let us at this point assume
  that the spectral triple
   is $QC^\infty$, i.e. ${\cla^\infty}$ and $\{ [D,a], ~a \in {\cla^\infty} \}$ are contained in the domains of all powers of the derivation $[|D|, \cdot]$. Under this assumption,
    $\tau$ is a positive faithful trace on the  $C^\ast$-subalgebra generated by ${\cla^\infty}$ and $\{ [D,a]~a \in {\cla^\infty} \}$,
     and the GNS  Hilbert space $L^2({\cla^\infty}, \tau)$ is denoted by $\clh^0_D$. Similarly, we equip $\Omega^1_D$ with a semi-inner product given by
     $<\eta,\eta^\prime>:=\tau(\eta^*\eta^\prime)$, and denote the Hilbert space obtained from it by $\clh^1_D$.
    The map $d_D : \clh^0_D \raro \clh^1_D$ given by $d_D(\cdot)=[D, \cdot]$ is an unbounded densely defined linear map.
    Let us  assume the following: \vspace{1mm}\\

\noindent{\bf Assumption}(i) (a) $d_D$ is closable
   (the
closure is denoted again by $d_D$);\\
(b) ${\cla^\infty} \subseteq {\rm Dom}(\cll) $, where
$\cll:=-d_D^*d_D$ and ${\cla^\infty}$ is viewed as a dense
subspace of $\clh^0_D$;\\

At this point, let us show that this assumption is valid under a
very natural condition on the spectral triple.   \vspace{1mm}\\
\blmma 
\label{d_closable}
Suppose that for every element $a \in {\cla^\infty}$, the map $\IR
\ni t \mapsto \al_t(X):={\rm exp}(itD) X {\rm exp}(-itD)$ is
differentiable at $t=0$ in the norm-topology of $\clb(\clh)$,
where $X=a$ or $[D,a]$. Then the assumption (i) is satisfied.
Moreover, in this case, $\cll$ maps ${\cla^\infty}$ into the weak closure of ${\cla^\infty}$ in $\clb(\clh^0_D)$. \elmma {\it Proof :}\\
We first observe that $\tau(\al_t(A))=\tau(A)$ for all $t$ and for
all $A \in \clb(\clh)$, since ${\rm exp}(itD)$ commutes with
$|{D}|^{-p}$.  If moreover, $A$  belongs to the domain of
norm-differentiability (at $t=0$) of $\al_t$, i.e.
$\frac{\al_t(A)-A}{t} \raro i[D,A]$ in operator-norm, then it
follows from the property of the Dixmier trace that
$\tau([D,A])=\frac{1}{i} \lim_{t \raro 0}
\frac{\tau(\al_t(A))-\tau(A)}{t}=0$. Now, since by assumption  we
have the norm- differentiability at $t=0$ of $\al_t(A)$ for $A $
belonging to the $\ast$-subalgebra (say $\clb$) generated by
${\cla^\infty}$ and $[D, {\cla^\infty}]$, it follows that $\tau([D,A])=0$ $\forall A
\in \clb$. Let us now fix $a,b,c \in {\cla^\infty}$ and observe that
\bean \lefteqn{<a~d_D(b),d_D(c)>}\\
&=& \tau((a~d_D(b))^*d_D(c)>\\
&=& -\tau([D,[D,b^*]a^*c])+\tau([D,[D,b^*]a^*]c)\\
&=& \tau([D,[D,b^*]a^*]c),\eean
using the fact that $\tau([D,[D,b^*]a^*c])=0$. This implies
$$|<a~d_D(b),d_D(c)>| \leq \| [D,[D,b^*]a^*] \|
\tau(c^*c)^{\frac{1}{2}}= \| [D,[D,b^*]a^*] \| \| c \|_2,$$ where
$\| c \|_2=\tau(c^*c)^{\frac{1}{2}}$ denotes the $L^2$-norm of $c
\in \clh^0_D$. This proves that $a~d_D(b)$ belongs to the domain
of $d_D^*$ for all $a,b \in {\cla^\infty}$, so in particular $d_D^*$ is
dense, i.e. $d_D$ is closable. Moreover, taking $a=1$, we see that
$d_D({\cla^\infty}) \subseteq {\rm Dom}(d_D^*)$, or in other words, ${\cla^\infty}
\subseteq {\rm Dom}(d_D^*d_D)$. This proves  (i)(a) and (i)(b).
The last sentence in the statement of  the lemma  can be proven
along the line of Theorem 2.9, page 129, \cite{fro}. \qed

We need few more assumptions on the operator $\cll$ to define the
quantum isometry
group.\vspace{1mm}\\
{\bf Assumption }(ii): $\cll$ has compact resolvents,\\
{\bf Assumption}(iii):  $\cll({\cla^\infty}) \subseteq {\cla^\infty}$;  \\
{\bf Assumption}(iv): Each eigenvector of $\cll$ (which has a discrete spectrum, hence a complete set of eigenvectors) belongs to ${\cla^\infty}$;\\
{\bf Assumption}(v)(`connectedness assumption'): the kernel of
$\cll$ is one-dimensional, spanned by the identity $1$ of ${\cla^\infty}$,
 viewed as a unit vector in $\clh^0_D$.  \vspace{2mm}\\
We call $\cll$ the noncommutative Laplacian and $T_t$ the
noncommutative heat semigroup. We summarize some simple
observations in form of the following \blmma \label{lem1}
(a) If the assumptions (i)-(v) are valid, then for $x \in {\cla^\infty}$, we have $\cll(x^*)=(\cll(x))^*$.\\
(b) If $T_t:={\rm exp}(t \cll)$ maps $\clh^0_D$ into  ${\cla^\infty}$ for all $t>0$, the the assumption (iv) is satisfied.
\elmma
{\it Proof :}\\
It follows by simple calculation using the facts that $\tau$ is a
trace and  $d_D(x^*)=-(d_D(x))^*$ that \bean \lefteqn{\tau
(\cll(x^*)^* y)}\\
&=& -\tau(d_D(x)d_D(y))=
-\tau(d_D(y)d_D(x))=\tau((d_D(y^*))^*d_D(x))\\
&=& <y^*,\cll(x)>=\tau(y\cll(x))= \tau(\cll(x)y), \eean for all $y
\in {\cla^\infty}$. By density of ${\cla^\infty}$ in $\clh^0_D$ (a) follows. To
prove (b), we note that if $x \in \clh^0_D$ is an eigenvector of
$\cll$, say $\cll(x)=\lambda x$ $(\lambda \in \IC)$, then we have
$T_t(x)=e^{\lambda t} x$, hence $x=e^{-\lambda t} T_t(x) \in {\cla^\infty}$.
 \qed \vspace{1mm}\\

 Since by assumption, $\cll$ has a countable set of
eigenvalues each with finite multiplicity, let us denote them by
$\lambda_0=0, \lambda_1, \lambda_2,...$ with  $V_0=\IC ~1,
V_1,V_2,...$ be corresponding eigenspaces (finite dimensional),
and  for each $i$, let $\{ e_{ij},j=1,...,d_i \}$ be an
orthonormal basis of $V_i$. By Assumption (iv), $V_i \subseteq
{\cla^\infty}$ for each $i$, $V_i$ is closed under $\ast$, and moreover,
$\{ e^\ast_{ij},j=1,...,d_i \}$ is also an orthonormal basis for
$V_i$, since $\tau(x^*y)=\tau(yx^*)$ for $x,y \in {\cla^\infty}$. We also
make the following \vspace{1mm}\\
 {\noindent}{\bf Assumption } (vi) The complex linear span of
$\{ e_{ij},i=0,1,...; j=1,...,d_i \}$, say ${\cla_0^\infty}$,
is norm-dense in ${\cla^\infty}$.\\

 \bdfn We say that a spectral triple satisfying the
assumptions (i)-(vi) admissible. \edfn \brmrk We have just seen
that  classical spectral triple $({\cla^\infty}=C^\infty(M), \clh,D)$,
where $M $ is compact connected spin manifold, $\clh$ is the $L^2$
space of square integrable spinors and $D$ is the Dirac operator,
is indeed admissible in our sense.
 Later on we shall discuss how we can weaken the connectedness assumption as well, thus accommodating a general classical (commutative) spectral
  triple in our set-up. Moreover, the standard examples of
  noncommutative spectral triples, e.g. those on $\cla_\theta$,
  quantum Heisenberg manifold etc., do belong to the admissible
  class.
  \ermrk

 \blmma \label{tracepres}  Let us assume that the spectral triple $({\cla^\infty},\clh,D)$ is
  admissible. Let  $\Psi :
{\cla_0^\infty} \raro {\cla_0^\infty}$ be a (norm-) bounded linear map, such that
$\Psi(1)=1$, and $\Psi \circ \cll=\cll \circ \Psi$ on the subspace
 ${\cla_0^\infty}$  spanned (algebraically) by $V_i,
i=1,2,...$. Then $\tau(\Psi(x))=\tau(x)$ for all $x \in {\cla^\infty}$.
\elmma
{\it Proof :}\\ By Lemma \ref{tracepres_1} with $\clh_1=\clh_2=\clh^0_D$, $\xi_1=\xi_2=1$, we have  $\tau(\Psi(x))=\tau(x)$ for all
$x \in {\cla_0^\infty}$. By the norm-continuity of $\Psi$ and $\tau$ it
extends to the whole of ${\cla^\infty}$. 
\qed \vspace{1mm}\\

\subsection{Definition and existence of the quantum isometry group}
 We begin by   recalling the definition of compact quantum groups and their actions from   \cite{woro}.  A 
compact quantum group is given by a pair $(\cls, \Delta)$, where $\cls$ is a unital separable $C^*$ algebra 
equipped
 with a unital $C^*$-homomorphism $\Delta : \cls \raro \cls \otimes \cls$ (where $\otimes$ denotes the injective tensor product)
  satisfying \\
  (ai) $(\Delta \ot id) \circ \Delta=(id \ot \Delta) \circ \Delta$ (co-associativity), and \\
  (aii) the linear span of $ \Delta(\cls)(\cls \ot 1)$ and $\Delta(\cls)(1 \ot \cls)$ are norm-dense in $\cls \ot \cls$. \\
  It is well-known (see \cite{woro}) that there is a canonical dense $\ast$-subalgebra $\cls_0$ of $\cls$, consisting of the matrix coefficients of
   the finite dimensional unitary (co)-representations of $\cls$, and maps $\epsilon : \cls_0 \raro \IC$ (co-unit) and
   $\kappa : \cls_0 \raro \cls_0$ (antipode)  defined
    on $\cls_0$ which make $\cls_0$ a Hopf $\ast$-algebra.

    We say that  the compact quantum group $(\cls,\Delta)$ acts on a unital $C^*$ algebra $\clb$,
    if there is a  unital $C^*$-homomorphism $\alpha : \clb \raro \clb \ot \cls$ satisfying the following :\\
    (bi) $(\alpha \ot id) \circ \alpha=(id \ot \Delta) \circ \alpha$, and \\
    (bii) the linear span of $\alpha(\clb)(1 \ot \cls)$ is norm-dense in $\clb \ot \cls$.\\
     \vspace{1mm}\\

 Let us now recall the concept of universal quantum groups as in
\cite{univ1}, \cite{free}
and references therein. We shall use most of the terminologies of
\cite{free}, e.g. Woronowicz $C^*$ -subalgebra, Woronowicz
$C^*$-ideal etc, however with the exception that we shall call the
Woronowicz $C^*$ algebras just compact quantum groups, and not use
the term compact quantum groups for the dual objects as done in
\cite{free}. For $Q \in GL_n(\IC)$, let $A_u(Q)$ denote the
universal compact quantum group generated by $u_{ij}, i,j=1,...,n$
satisfying the relations
$$u u^*=I_n =u^*u, ~~u^\prime Q \overline{u} Q^{-1}=I_n=Q
\overline{u} Q^{-1} u^\prime,$$ where $u=(( u_{ij} )), $
$u^\prime=(( u_{ji} ))$ and $\overline{u}=(( u_{ij}^* ))$. The coproduct, say $\tilde{\Delta}$, is given by, $$ \tilde{\Delta}(u_{ij})=\sum_k u_{ik} \ot u_{kj}.$$ We
refer the reader to \cite{univ1}  for a detailed  discussion on the structure and classification of
such quantum groups. Let us denote by $\clu_i$ the quantum group
$A_{d_i}(I)$, where $d_i$ is dimension of the subspace $V_i$. We
fix a representation $\beta_i : V_i \raro V_i \otimes \clu_i$ of
$\clu_i$ on the Hilbert space $V_i$, given by $\beta_i ( e_{ij})=
\sum_k e_{ik} \otimes u^{(i)}_{kj}$, for $j=1,...,d_i$, where
$u^{(i)}\equiv u^{(i)}_{kj}$ are the generators of $\clu_i$ as
discussed before. Thus, both $u^{(i)}$ and $\bar{u^{(i)}}$ are
unitaries.  It follows from \cite{free} that the representations
$\beta_i$ canonically induce a representation $\beta=\ast_i
\beta_i$ of the free product $\clu:=\ast_i \clu_i$ (which is a
compact quantum group, see \cite{free} for the details) on the
Hilbert space $\clh^0_D$, such that the restriction of $\beta$ on
$V_i$ coincides with $\beta_i$ for all $i$.

In view of the  characterization of smooth isometric action on  a classical manifold, we make the following definitions.
         \bdfn A quantum family of smooth isometries of a 
noncommutative manifold ${\cla^\infty}$ (or, more precisely on the
corresponding  spectral triple) is a pair $(\cls, \alpha)$ where $\cls$ is a separable unital $C^*$-algebra,  $\alpha : \overline{\cla} \raro
\overline{\cla} \otimes \cls$ (where $\overline{\cla}$ denotes the
$C^*$ algebra obtained by completing ${\cla^\infty}$ in the norm of
$\clb(\clh^0_D)$) is a unital $C^*$-homomorphism, satisfying the following:\\
(a) $ \overline{{\rm Sp}}\left( \alpha(\overline{\cla})(1 \ot \cls) \right)=\overline{\cla} \ot \cls$,\\
(b)  $\alpha_\phi:=(id \otimes \phi)\circ
\alpha$ maps ${\cla^\infty_0}$ into itself and commutes with $\cll$ on
$\cla^\infty_0$, for every state $\phi$ on $\cls$.\\ 

In case the $C^*$-algebra $\cls$ has a coproduct $\Delta$ such that $(\cls,\Delta)$ is a compact quantum group and $\alpha$ is an action of $(\cls, \Delta)$ on $\overline{\cla}$, we say that $(\cls, \Delta)$ acts smoothly and isometrically on the noncommutative manifold.
\edfn

Fix a  spectral triple $(\cla^\infty, \clh, D)$. Consider the category ${\bf Q}$ with the object-class consisting of all quantum families of isometries $(\cls, \alpha)$ of the given noncommutative manifold, and the set of morphisms ${\rm Mor}((\cls,\alpha),(\cls^\prime,\alpha^\prime))$ being the set of unital $C^*$-homomorphisms $\phi : \cls \raro \cls^\prime$ satisfying $({\rm id} \ot \phi) \circ \alpha=\alpha^\prime$. We also consider another category ${\bf Q}^\prime$ whose objects are triplets $(\cls, \Delta, \alpha)$ where $(\cls,\Delta)$ is a compact quantum group acting smoothly and isometrically on the given noncommutative manifold, with $\alpha$ being the corresponding action. The morphisms  are the homomorphisms of compact quantum groups which are also morphisms of the underlying quantum families. The forgetful functor $F: {\bf Q}^\prime \raro {\bf Q}$ is clearly faithful, and we can view $F({\bf Q}^\prime)$ as a subcategory of ${\bf Q}$. 

Let us assume from now on that the spectral triple $(\cla^\infty, \clh, D)$ is admissible.  Our aim is to prove the existence of a universal object in ${\bf Q}$. We shall also  prove that the (unique upto isomorphism) universal object  belongs to $F({\bf Q}^\prime)$, and its pre-image in ${\bf Q}^\prime$ is a universal object in the category ${\bf Q}^\prime$. To this end, we need some preparatory results.
 \blmma
\label{lem2} Consider an admissible spectral triple $(\cla^\infty,\clh,D)$ and let $(\cls,\alpha)$ be a quantum family of smooth isometries of the spectral triple.  Moreover, assume that the
action  $\alpha$ is faithful in the sense that there is no
proper  $C^*$-subalgebra $\cls_1$ of $\cls$ such that
$\alpha({\cla^\infty}) \subseteq {\cla^\infty} \otimes \cls_1$. 
Then $\tilde{\alpha}: \cla^\infty \ot \cls \raro \cla^\infty \ot \cls$ defined by $\tilde{\alpha}(a \ot b): \alpha(a)(1 \ot b)$ extends to an $\cls$-linear unitary on the Hilbert $\cls$-module $\clh^0_D \ot \cls$, denoted again by $\tilde{\alpha}$. 
 Moreover, we
can find a $C^*$-isomorphism  $\phi : \clu/
\cli \raro \cls$ between $\cls$ and a quotient of $\clu$ by a
 $C^*$-ideal
 $\cli$ of $\clu$, such that $ \alpha= ({\rm id}\ot \phi) \circ ({\rm id} \ot \Pi_\cli) \circ
 \beta$ on ${\cla^\infty} \subseteq \clh^0_D$, where $\Pi_\cli$ denotes the quotient map from $\clu$ to
 $\clu/\cli$.

 If, furthermore, there is a compact quantum group structure on $\cls$ given by a coproduct  $\Delta$ such that $(\cls,\Delta, \alpha)$ is an object in ${\bf Q}^\prime$, the map  $\alpha : {\cla^\infty}
\raro {\cla^\infty} \otimes \cls$ extends to a unitary representation
(denoted again by $\alpha$) of the compact quantum group $(\cls,\Delta)$ on $\clh^0_D$. In this case, the ideal $\cli$ is a Woronowicz $C^*$-ideal and the $C^*$-isomorphism $\phi : \clu/ \cli \raro \cls$ is a morphism of compact quantum groups.
\elmma
{\it Proof :}\\
Let $\omega$ be any state on $\cls$. Since the action $\alpha :
{\cla^\infty} \raro {\cla^\infty} \otimes \cls$ is smooth and isometric, we conclude
by Lemma \ref{tracepres} that $\tau(\alpha_\omega(x))=\tau(x)
\omega(1)$ for all $x \in \overline{\cla}$. Since $\omega$ is
arbitrary, we have $(\tau \otimes id) \alpha(x)=\tau(x) 1_\cls$
for all $x \in \overline{\cla}$. So, $<\alpha(x),
\alpha(y)>_{\cls}=<x,y> 1_\cls$, where $<\cdot,\cdot>_{ \cls}$
denotes the $\cls$-valued inner product of the Hilbert module
$\clh^0_D \ot \cls$. This proves that $\tilde{\alpha} $ defined by
$\tilde{\alpha}(x \ot b):=\alpha(x)(1 \ot b)$ ($x \in {\cla^\infty}, b \in
\cls$) extends to an $\cls$-linear isometry on the Hilbert
$\cls$-module $\clh^0_D \otimes \cls$. Moreover, since
$\alpha({\cla^\infty})(1 \ot \cls)$ is norm-dense in $\bar{\cla} \ot \cls$,
it is clear that the $\cls$-linear span of the range of
$\alpha({\cla^\infty})$ is dense in the Hilbert module $\clh^0_D \ot \cls$,
or in other words, the isometry $\tilde{\alpha}$ has a dense
range, so it is a unitary. 

Since $\alpha_\omega$ leaves each $V_i$
invariant, it is clear that $\alpha$  maps $V_i$ into $V_i \otimes
\cls$ for each $i$. Let $v^{(i)}_{kj}$ ($j,k=1,...,d_i$) be the
elements of $\cls$ such that $\alpha(e_{ij})=\sum_k e_{ik} \otimes
v^{(i)}_{kj}$.  Note that $v_i:=((v^{(i)}_{kj} ))$ is a unitary in
$M_{d_i}(\IC) \otimes \cls$. Moreover, the $\ast$-subalgebra
generated by all $ \{ v^{(i)}_{kj}, i,j,k \geq 1 \}$ must be
dense in $\cls$ by the assumption of faithfulness.

 We have already remarked that $\{ e_{ij}^* \}$ is also an
orthonormal basis of $V_i$, and since $\alpha$, being a
$C^*$-action on $\overline{\cla}$, is $\ast$-preserving, we have
$\alpha(e^*_{ij})=(\alpha(e_{ij}))^*=\sum_k e^*_{ik} \otimes
v^{(i)^*}_{kj}$, and therefore $(( v^{(i)^*}_{kj}))$ is also
unitary. By universality of $\clu_i$, there is a
$C^*$-homomorphism from $\clu_i$ to $\cls$ sending $u^{(i)}_{kj} $
to $v^{(i)}_{kj}$, and by definition of the free product, this
induces a $C^*$-homomorphism, say $\Pi$, from $\clu$ onto $\cls$,
 so that $\clu/\cli
\cong \cls$, where  $\cli:={\rm Ker}(\Pi)$.

In case $\cls$ has a coproduct $\Delta$ making it into a compact quantum group and $\alpha$  is a quantum group action, it is easy to see that  the subalgebra of $\cls$ generated by $v^{(i)}_{kj}$ is a Hopf
algebra, with $\Delta(v^{(i)}_{kj})=\sum_l v^{(i)}_{kl} \ot v^{(i)}_{lj}$. From this, it follows that  $\Pi$ is Hopf-algebra morphism, hence $\cli$ is a Woronowicz $C^*$-ideal.
\qed \vspace{1mm}\\

Before we state and prove the main theorem, let us note the following elementary fact about $C^*$-algebras.

  \blmma \label{lim} Let $\clc$ be a
$C^*$ algebra and $\clf$ be a nonempty collection of $C^*$-ideals (closed two-sided ideals) of
$\clc$. Then for any $x \in \clc$, we have
$$ \sup_{I \in \clf}   \| x+I \| = \|x+ I_0 \|,$$
where $I_0$ denotes the intersection of all $I$ in $\clf$ and $\|
x+I\|=inf \{ \| x -y \|~:~y \in I \}$ denotes the norm in
$\clc/I$. \elmma
 {\it Proof :}\\
It is clear that $\sup_{I \in \clf}   \| x+I \|$ defines a norm on
$\clc/I_0$, which is in fact a $C^*$-norm since each of the
quotient norms  $\|\cdot+I\|$ is so. Thus the lemma follows from
the uniqueness of $C^*$ norm on the $C^*$ algebra $\clc/I_0$. \qed
 \bthm
\label{main} For any admissible spectral triple $(\cla^\infty, \clh, D)$, the category ${\bf Q}$ of quantum families of smooth isometries has a universal (initial) object, say  $(\clg, \alpha_0)$. Moreover, $\clg$ has a coproduct $\Delta_0$ such that $(\clg,\Delta_0)$ is a compact quantum group and $(\clg,\Delta_0,\alpha_0)$ is a universal object in the category ${\bf Q}^\prime$ of compact quantum groups acting smoothly and isometrically on the given spectral triple.  The action $\alpha_0$ is faithful. \ethm
 {\it Proof :}\\
Recall the $C^*$-algebra $\clu$ considered before, and the map  $\beta$ 
from  $ \clh^0_D$ to $\clh_D^0 \ot \clu$.   By our
definition of $\beta$, it is clear that $\beta({\cla_0^\infty}) \subseteq
{\cla_0^\infty} \otimes_{\rm alg} \clu$. However, $\beta$ is only a linear
map (unitary) but not necessarily a $\ast$-homomorphism. We shall
construct the universal object as a suitable quotient of $\clu$.
Let $\clf$ be the collection of all those  $C^*$-ideals
$\cli$ of $\clu$ such that the composition $\Gamma_\cli:=(id
\otimes \Pi_\cli) \circ \beta : {\cla_0^\infty} \raro {\cla_0^\infty} \otimes_{\rm
alg} (\clu/\cli)$ extends to a $C^*$-homomorphsim from
$\bar{\cla}$ to $\bar{\cla} \ot (\clu/\cli)$, where $\Pi_\cli$
denotes the quotient map from $\clu$ onto $\clu/\cli$. This
collection is nonempty, since the trivial one-dimensional $C^*$-algebra $\IC$ gives an object in ${\bf Q}$ and by Lemma \ref{lem2} we do get a member of $\clf$. Now,
let $\cli_0$ be the intersection of all ideals in $\clf$. We claim
that $\cli_0$ is again a member of $\clf$.   Since any
$C^*$-homomorphism is contractive, we have $\| \Gamma_\cli(a) \|
\equiv \|\beta(a)+ \bar{\cla} \ot \cli\| \leq \| a \|$ for all $a
\in {\cla_0^\infty}$ and $\cli \in \clf$. By Lemma \ref{lim}, we see that
$\|\Gamma_{\cli_0}(a) \| \leq \| a \|$ for $a \in {\cla_0^\infty}$, so
$\Gamma_{\cli_0}$ extends to a norm-contractive map on
$\bar{\cla}$ by the density of ${\cla_0^\infty}$ in $\bar{\cla}$. Moreover, 
for $a,b \in \bar{\cla}$ and for $\cli \in \clf$, we have
$\Gamma_{\cli}(ab)=\Gamma_{\cli}(a)\Gamma_{\cli}(b)$. Since
$\Pi_{\cli}=\Pi_{\cli} \circ \Pi_{\cli_0}$, we can rewrite the
homomorphic property of $\Gamma_{\cli}$ as
$$\Gamma_{\cli_0}(ab)-\Gamma_{\cli_0}(a) \Gamma_{\cli_0}(b) \in \bar{\cla} \otimes (\cli/\cli_0).$$
Since this holds for every $\cli \in \clf$, we conclude that
$\Gamma_{\cli_0}(ab)-\Gamma_{\cli_0}(a) \Gamma_{\cli_0}(b) \in
\bigcap_{\cli \in \clf} \bar{\cla} \otimes (\cli/\cli_0)=(0)$,
i.e. $\Gamma_{\cli_0}$ is a homomorphism.  In a similar way, we
can show that it is a $\ast$-homomorphism. 
Since each $\beta_i$ is a unitary representation of the compact quantum group $\clu_i$ on the finite dimensional space $V_i$, it follows that $\beta_i(V_i)(1 \ot \clu_i)$ is total in $V_i \ot \clu_i$. In particular, for any $v_i \in V_i$ ($i$ arbitrary), the element $v_i \ot 1_{\clu_i}=v_i \ot 1_{\clu}$ belongs to the linear span of $\beta_i(V_i)(1 \ot \clu_i) \subset \beta(V_i)(1 \ot \clu)$. Thus, $\cla_0^\infty \ot 1_\clu$  is contained in the  linear span of $\beta(\cla^\infty_0)(1 \ot \clu)$ and hence $\cla^\infty_0 \ot 1_{\frac{\clu}{\cli_0}}$ is linearly spanned by $\Gamma_{\cli_0}(\cla^\infty_0)(1 \ot \clu/\cli_0)$.  By the norm-denisty of $\cla^\infty_0$ in $\overline{\cla}$ and the contractivity of the quotient map, it follows that $\overline{\cla} \ot \clu/\cli_0$ is the closed linear span of $\Gamma_{\cli_0}(\cla^\infty_0)(1 \ot \clu/\cli_0)$. This completes the proof that $(\clu/\cli_0, \Gamma_{\cli_0})$ is indeed an object of ${\bf Q}$.

We now show   that $\clg:=\clu/\cli_0$ is a 
universal object in ${\bf Q}$. To see this, consider any object $(\cls,\alpha)$ of ${\bf Q}$.  Without loss
of generality we can assume the action to be faithful, since
otherwise we can replace $\cls$ by the  $C^*$-subalgebra
generated by the elements $\{ v^{(i)}_{kj} \}$ appearing in the proof of Lemma \ref{lem2}.  But by
Lemma \ref{lem2} we can further assume that $\cls$ is isomorphic
with $\clu/\cli$ for some $\cli \in \clf$. Since $\cli_0 \subseteq
\cli$, we have a $C^*$-homomorphism from
$\clu/\cli_0$ onto $\clu/\cli$, sending $x+\cli_0$ to $x +\cli$, which is clearly a morphism in the category ${\bf Q}$. This is indeed the unique such morphism, since it is uniquely determined on the dense subalgebra generated by $\{ u^{(i)}_{kj}+\cli_0,~i,j,k \geq 1 \}$ of $\clg$.

To construct the coproduct on $\clg=\clu/\cli_0$, we first consider $\alpha^{(2)}=(\Gamma_{\cli_0} \ot {\rm id}) \circ \Gamma_{\cli_0} : \overline{\cla} \raro \overline{\cla} \ot \clg \ot \clg$. It is easy to verify that $(\clg \ot \clg, \alpha^{(2)})$ is an object in the category ${\bf Q}$, so by the universality of $(\clg,\Gamma_{\cli_0})$, we have a unique unital $C^*$-homomorphism $\Delta_0 : \clg \raro \clg \ot \clg$ satisfying $$ ({\rm id} \ot \Delta_0) \circ \Gamma_{\cli_0}(x) =\alpha^{(2)}(x)~\forall x \in \overline{\cla}.$$ 
Taking $x=e_{ij}$, we get $$ \sum_l e_{il} \ot (\pi_{\cli_0} \ot \pi_{\cli_0}) \left( \sum_k u^{(i)}_{lk} \ot u^{(i)}_{kj} \right)=\sum_l e_{il} \ot \Delta_0(\pi_{\cli_0}(u^{(i)}_{lj})).$$ Comparing coefficients of $e_{il}$, and recalling that $\tilde{\Delta}(u^{(i)}_{lj})=\sum_k u^{(i)}_{lk} \ot u^{(i)}_{kj}$ (where $\tilde{\Delta}$ denotes the coproduct on $\clu$), we have \be \label{coprod}(\pi_{\cli_0} \ot \pi_{\cli_0}) \circ \tilde{\Delta}=\Delta_0 \circ \pi_{\cli_0} \ee on the linear span of $\{ u^{(i)}_{jk}, i,j,k \geq 1 \}$, and hence on the whole of $\clu$. 
This implies that $\Delta_0$ maps $\cli_0={\rm Ker}(\pi_{\cli_0})$ into ${\rm Ker}(\pi_{\cli_0} \ot \pi_{\cli_0})=(\cli_0 \ot 1 +1 \ot \cli_0) \subset \clu \ot \clu$. In other words, $\cli_0$ is a Hopf $C^*$-ideal, and hence $\clg=\clu/\cli_0$ has the canonical compact quantum group structure as a quantum subgroup of $\clu$. It is clear from the relation (\ref{coprod}) that  $\Delta_0$ coincides with the canonical coproduct of the quantum subgroup $\clu/\cli_0$ inherited from that of $\clu$. It is also easy to see that the object $(\clg, \Delta_0, \Gamma_{\cli_0})$ is universal in the category ${\bf Q}^\prime$, using the fact that (by Lemma \ref{lem2}) any compact quantum group $(\clg, \Phi)$ acting smoothly and isometrically on the given spectral triple is isomorphic with a quantum subgroup $\clu/\cli$, for some Hopf $C^*$-ideal $\cli$ of $\clu$. 

Finally, the faithfulness of $\alpha_0$ follows from the universality by standard arguments which we briefly sketch. If $\clg_1 \subset \clg$ is a $\ast$-subalgebra of $\clg$ such that $\alpha_0(\overline{\cla}) \subseteq \overline{\cla} \ot \clg_1$, it is easy to see that $(\clg_1, \Delta_0, \alpha_0)$ is also a universal object, and by definition of universality of $\clg$ it follows that there is a unique morphism, say $j$, from $\clg$ to $\clg_1$. But the map $j \circ i$ is a morphism from $\clg$ to itself, where $i : \clg_1 \raro \clg$ is the inclusion. Again by universality, we have that $j \circ i ={\rm id}_\clg$, so in particular, $i$ is onto, i.e. $\clg_1=\clg$.  
\qed

\bdfn We shall call the universal object $(\clg,\Delta_0)$
obtained in the theorem above the quantum isometry group of
$({\cla^\infty},\clh,D)$ and denote it by $QISO({\cla^\infty},\clh,D)$, or just
$QISO({\cla^\infty})$ (or sometimes $QISO(\bar{\cla})$)  if the spectral
triple is understood from the context. \edfn
\brmrk
\label{a_0^inf}
Assume that an admissible spectral triple $(\cla^\infty, \clh,D)$ also satisfies the condition (i) of Lemma \ref{classic2}, i.e. $\bigcap {\rm Dom}(\cll^n)=\cla^\infty$. Let $\alpha : \overline{\cla} \raro \overline{\cla} \ot \cls$ be a smooth isometric action on $\cla^\infty$ by a compact quantum group $\cls$. We recall from the proof of Lemma \ref{lem2} that the map $\tilde{\alpha}$ from $ \overline{\cla} \ot \cls$ to itself extends to an $\cls$-linear  unitary on the Hilbert $\cls$-module $\clh^0_D \ot \cls$, i.e.  $\tilde{\alpha}$ can be viewed as a unitary in $\clb(\clh^0_D) \ot \cls$. Clearly, for any state $\phi$ on $\cls$, we have $\alpha_\phi=({\rm id} \ot \phi)(\tilde{\alpha}) \in \clb(\clh^0_D)$. Now, by the definition of a smooth isometric action, the bounded operator $\alpha_\phi$ commutes with  the self-adjoint operator $\cll$ on $\cla^\infty_0$, which is a core for $\cll$. So, $\alpha_\phi$ must commute with $\cll^n$ for all $n$, and in particular keeps $\cla^\infty=\bigcap_n {\rm Dom}(\cll^n)$ invariant.
\ermrk

\brmrk
Let us now briefly indicate how one can weaken the hypothesis of
connectedness. Such an extension of our results is desirable to
accommodate  the classical spaces, including the finite sets and
graphs, in our framework.  
 One possibile approach  could be to consider the category of compact quantum group actions  $\alpha$ which are not only `smmoth' and `isometric' in our sense, but  also satisfy the $\tau$-invariance condition, i.e. $(\tau \ot {\rm id})(\alpha(a))=\tau(a) 1$. It is easy to see that the connectedness assumption has been used by us only to prove that the $\tau$-invariance is automatic for smooth isometric actions. Thus, if we work in the smaller cateogory of such $\tau$-invraiant actions only, the proof of Theorem \ref{main} does go through and thus we can prove the existence of a universal object, to be defined as the quantum isometry group.  
It is easy to see that for the algebra of functions on a finite set, with the spectral triple given by $D=0$, this quantum isometry group  coincides with thw quantum permutation group defined by Wang. 
\ermrk

\brmrk It is easy to see how to extend our formulation and results
to spectral triples which are not necessarily of type II, i.e.
when the trace $\tau$ is replaced by some non-tracial positive
functional. Indeed, our construction will go through in such a
situation more or less verbatim, by replacing the universal
quantum groups $A_{d_i}(I)$ by $A_{d_i}(Q_i)$ for some suitable
choice of matrices $Q_i$ coming from the modularity property of
$\tau$. \ermrk

\subsection{Construction of quantum group-equivariant spectral triples}
In this subsection, we shall briefly discuss the relevance of quantum isometry group to the problem of constructing quantum group equivariant spectral triples, which is important to understand the role of quantum groups in the framework of noncommutative geometry. There has been a lot of activity in this direction recently, see, for example, the articles by Chakraborty and Pal (\cite{partha}),  Connes (\cite{con2}), Landi et al (\cite{landiqgp}) and the references therein. 
In the classical situation, there exists a natural unitary representation of the isometry group $G=ISO(M)$  of a manifold $M$ on the Hilbert space of forms, so that the operator $d+d^*$ (where $d$ is the de-Rham differential operator) commutes with the representation. Indeed, $d+d^*$ is also a Dirac operator for the spectral triple given by the natural representation of $C^\infty(M)$ on the Hilbert space of forms, so we have a canonical construction of $G$-equivariant spectral triple. Our aim in this subsection is to generalize this to the noncommutative framework, by proving that $d_D+d_D^*$ is equivariant with respect to a canonical unitary representation  on the Hilbert space of `noncommutative forms' (see, for example, \cite{fro} for a detailed discussion of such forms). 

Consider an admissible spectral triple $(\cla^\infty, \clh, D)$ and moreover, make the assumption of Lemma \ref{d_closable}, i.e. assume that $t \mapsto e^{itD} x e^{-itD}$ is norm-differentiable at $t=0$ for all $x$ in the $\ast$-algebra $\clb$ generated by $\cla^\infty$ and $[D, \cla^\infty]$.  

\blmma
\label{d_*}
In the notation of Lemma \ref{d_closable}, we have the following (where $b,c \in \cla^\infty$): \be \label{***} d_D^*(d_D(b)c)=-\frac{1}{2} \left( b\cll(c)-\cll(b)c-\cll(bc) \right).\ee
\elmma
{\it Proof:}\\
Denote by $\chi(b,c)$ the right hand side of euqation (\ref{***}) and fix any $a \in \cla^\infty$. Using the facts the the functional $\tau$ is a faithful trace on the $\ast$-algebra $\clb$,  $\cll=-d_D^*d_D$ and that  $[D, X]=0$ for any $X$ in $\clb$, we have, \bean \lefteqn{\tau(a^* \chi(b,c))}\\
&=& -\frac{1}{2} \{ \tau(a^*b \cll(c))+\tau(ca^* \cll(b))+\tau(a^* \cll(bc))\}\\
&=& \frac{1}{2} \{ \tau([D,a^*b][D,c])-\tau([D,ca^*][D,b])-\tau([D,a^*][D,bc])\}\\
&=& \frac{1}{2} \{ \tau(a^*[D,b][D,c])-\tau([D,c]a^*[D,b])-\tau(c[D,a^*][D,b])-\tau([D,a^*][D,b]c)\}\\
&=& -\tau([D,a^*][D,b]c)\\
&=&\tau([D,a]^*[D,b]c)\\
&=&\lgl d_D(a), d_D(b)c \rgl\\
&=& \tau(a^*(d_D^*(d_D(b)c))).
\eean
\qed\\

From this, we get the following by a simple computation: \be \label{innerprodformula}
\lgl a d_D(b), a^\prime d_D(b^\prime) \rgl=-\frac{1}{2} \tau(b^* \Psi(a^*a, b^\prime)),\ee for $a,b,a^\prime,b^\prime \in \cla^\infty$, and where $\Psi(x,y):=\cll(xy)-\cll(x)y+x\cll(y).$ Now, let us denote the quantum isometry group of the given spectral triple $(\cla^\infty, \clh, D)$ by $(\clg, \Delta, \alpha)$. Let $\cla_0$ denote the $\ast$-algebra generated by $\cla^\infty_0$, $\clg_0$ denote $\ast$-algebra of $\clg$ generated by matrix elements of irreducible representations. Clearly, $\alpha : \cla_0 \raro \cla_0 \ot_{\rm alg} \clg_0$ is a Hopf-algebraic action of $\clg_0$ on $\cla_0$. Define $\tilde{\Psi} : (\cla_0 \ot_{\rm alg} \clg_0) \times (\cla_0 \ot_{\rm alg} \clg_0) \raro \cla_0 \ot _{\rm alg} \clg_0$ by $$ \tilde{\Psi}((x \ot q),  (x^\prime \ot q^\prime)):=\Psi(x,x^\prime) \ot (qq^\prime).$$ It follows from the relation $(\cll \ot {\rm id}) \circ \alpha=\alpha \circ \cll$ on $\cla_0$ that \be \label{psi_L} \tilde{\Psi}(\alpha(x),\alpha(y))=\alpha(\Psi(x,y)).\ee We now define a linear map $\alpha^{(1)}$ from the linear span of $\{ a d_D(b):~a,b \in \cla_0 \}$ to $\clh^1_D \ot  \clg$ by setting $$ \alpha^{(1)}(ad_D(b)):=\sum_{i,j} a^{(1)}_id_D(b^{(1)}_j) \ot a^{(2)}_ib^{(2)}_j,$$ where for any $x \in \cla_0$ we write $\alpha(x)=\sum_i x^{(1)}_i \ot x^{(2)}_i \in \cla_0 \ot_{\rm alg} \clg_0$ (summation over finitely many terms). We shall sometimes use the Sweedler convention of writing the above simply as $\alpha(x)=x^{(1)} \ot x^{(2)}$.
It then follows from the identities (\ref{innerprodformula}) and (\ref{psi_L}), and also the fact that $(\tau \ot {\rm id}) ( \alpha(a))=\tau(a)1$ for all $a \in \cla_0$  that \bean \lefteqn{\lgl ad_D(b), a^\prime d_D(b^\prime) \rgl_\clg }\\
&=& -\frac{1}{2} (\tau \ot {\rm id})(\alpha(b^*)\tilde{\Psi}(\alpha(a^* a^\prime), \alpha(b^\prime)))\\
&=& -\frac{1}{2} (\tau \ot {\rm id})(\alpha(b^*) \alpha(\Psi(a^*a^\prime,b^\prime)))\\
&=& -\frac{1}{2} (\tau \ot {\rm id})(\alpha(b^* \Psi(a^*a^\prime, b^\prime)))\\
&=& -\frac{1}{2} \tau(b^* \Psi(a^*a^\prime, b^\prime)) 1_\clg \\
&=& \lgl ad_D(b), a^\prime d_D(b^\prime) \rgl 1_\clg. \eean
This proves that $\alpha^{(1)}$ is indeed well-defined and extends to a $\clg$-linear isometry on $\clh^1_D \ot \clg$, to be denoted by $U^{(1)}$, which sends $(a d_D(b)) \ot q$ to $\alpha^{(1)}(a d_D(b))(1 \ot q)$, $a,b \in \cla_0$, $q \in \clg$. Moreover, since the linear span of $\alpha(\cla^\infty_0)(1 \ot \clg)$ is dense in $\clh^0_D \ot \clg$, it is easily seen that the range of the isometry $U^{(1)}$ is the whole of $\clh^1_D \ot \clg$, i.e. $U^{(1)}$ is a unitary. In fact, from its definition it can also be shwon that $U^{(1)}$ is a unitary representation of the compact quantum group $\clg$ on $\clh^1_D$. 

In a similar way, we can construct unitary representation $U^{(n)}$ of $\clg$ on the Hilbert space of $n$-forms for any $n \geq 1$, by defining $$ U^{(n)}((a_0 d_D(a_1) d_D(a_2)...d_D(a_n)) \ot q)=a_0^{(1)}d_D(a_1^{(1)})...d_D(a_n^{(n)}) \ot (a_0^{(2)}a_1^{(2)}...a_n^{(2)}q),~~   a_i \in \cla_0^\infty,$$ (using Sweedler convention) and verifying that it extends to a unitary. We also denote by $U^{(0)}$ the unitary representation $\tilde{\alpha}$ on $\clh^0_D$ discussed before. Finally, we have a unitary representation $U=\bigoplus_{n \geq 0} U^{(n)}$ of $\clg$ on $\tilde{\clh}:=\bigoplus_n \clh^n_D$, and also extend $d_D$ as a closed densely defined operator on $\tilde{\clh}$ in the obvious way, by defining $d_D(a_0d_D(a_1)...d_D(a_n))=d_D(a_0)...d_D(a_n)$. It is now straightforward to see the following:\\
\bthm
The operator $D^\prime:=d_D+d_D^*$ is equivariant in the sense that $U (D^\prime \ot 1)=(D^\prime \ot 1)U$. 
\ethm
We point out that there is a natural representation $\pi$ of $\overline{\cla}$ on $\tilde{\clh}$ given by $\pi(a)(a_0d_D(a_1)...d_D(a_n))=aa_0d_D(a_1)...d_D(a_n)$, and $(\pi(\cla^\infty), \tilde{\clh}, D^\prime)$ is indeed a spectral triple, which is $\clg$-equivariant. 

Although the relation between spectral properties of $D$ and $D^\prime$ 
   is not clear in general, in many cases of interest  (e.g. when there is an underlying type $(1,1)$ spectral data in the sense of \cite{fro}) these two Dirac operators are closely related.  
      As an illustration, consider  the canonical spectral  on the noncommutative $2$-torus $\cla_\theta$, which is discussed in some details in the next section. In this case, the Dirac operator $D$ acts on $L^2(\cla_\theta, \tau) \ot \IC^2$, and it can easily be shown (see \cite{fro})  that the Hilbert space of forms is isomorphic with $L^2(\cla_\theta,\tau) \ot \IC^4 \cong L^2(\cla_\theta) \ot \IC^2 $; thus $D^\prime$ is essentially same as $D$ in this case.

\section{Examples and computations}
We give some simple yet interesting explicit examples of quantum
isometry groups here. However, we give only some  computational
details for the first example, and for the rest, the reader is
referred to   a companion article (\cite{jyotish}).\\
\vspace{1mm}\\
\noindent{\bf Example 1 : commutative tori}\\
Consider $M=\IT$, the one-torus, with the usual Riemannian
structure. The $\ast$-algebra ${\cla^\infty}=C^\infty(M)$ is generated by
one unitary $U$, which is the multiplication operator by $z$ in
$L^2(\IT)$.  The Laplacian is given by $\cll(U^n)=- n^2
U^n$. If a compact quantum group $(\cls, \Delta_\cls)$ acts on
${\cla^\infty}$ smoothly, let $A_n, n \in \IZ$ be elements of $\cls$ such
that $\alpha_0(U)=\sum_n U^n \ot A_n$ (here $\alpha_0 : {\cla^\infty} \raro
{\cla^\infty} \ot_{\rm alg} \cls$ is the $\cls$-action on ${\cla^\infty}$). Note
that this infinite sum converges at least in the topology of the
Hilbert space $L^2(\IT) \ot L^2(\cls)$, where $L^2(\cls)$ denotes
the GNS space for the Haar state of $\cls$.  It is  clear that the
 condition  $(\cll \ot id) \circ \alpha_0=\alpha_0 \circ \cll$ forces to have  $A_n=0$ for all but $n=\pm 1$. The
conditions $\alpha_0(U)\alpha_0(U)^*=\alpha_0(U)^* \alpha_0(U)=1
\ot 1$ further imply  the following:
$$A_1^*A_1+A_{-1}^*A_{-1}=1=A_1A_1^*+A_{-1}A_{-1}^*,$$ $$
A_1^*A_{-1}=A_{-1}^*A_1=A_1A_{-1}^*=A_{-1}A_1^*=0.$$ It follows
that $A_{\pm 1}$ are partial isometries with orthogonal domains
and ranges. Say, $A_1$ has domain $P$ and range $Q$. Hence the
domain and range of $A_{-1}$ are respectively $1-P$ and $1-Q$.
Consider the unitary $V=A+B$, so that $VP=A$, $V(1-P)=B$. Now,
from the fact that $(\cll \ot id)
(\alpha_0(U^2))=\alpha_0(\cll(U^2))$ it is easy to see that the
coefficient of $1 \ot 1$ in the expression of $\alpha_0(U)^2$ must
be $0$, i.e. $AB+BA=0$. From this, it follows that $V$ and $P$
commute and therefore $P=Q$. By straightforward calculation  using
the facts that $V$ is unitary, $P$ is a projection and $V$ and $P$
commute, we can verify that $\alpha_0 $ given by $\alpha_0(U)=U
\ot VP+U^{-1}\ot V(1-P)$ extends to a $\ast$-homomorphsim from
${\cla^\infty}$ to ${\cla^\infty} \ot C^*(V,P)$ satisfying $(\cll \ot id) \circ
\alpha_0=\alpha_0 \circ \cll$. It follows that the $C^*$ algebra
$QISO(\IT)$ is commutative and generated by a unitary $V$ and a
projection $P$, or equivalently by
 two partial isometries $A$, $B$ such that $A^*A=AA^*, B^*B=BB^*, AB=BA=0$.
 So, as a $C^*$   algebra it is isomorphic with $C(\IT) \oplus C(\IT) \cong C(\IT \times \IZ_2)$.
 The coproduct (say $\Delta_0$) can easily be calculated from the
 requirement of co-associativity, and
   the Hopf algebra structure of  $QISO(\IT)$ can be seen to coincide with that of
  the semi-direct product of $\IT$ by $\IZ_2$, where the generator of $\IZ_2$ acts on $\IT$ by sending $z \mapsto \bar{z}$.

  We summarize this in form of the following.
\bthm The universal quantum group of isometries $QISO(\IT)$ of the
one-torus $\IT$ is isomorphic (as a quantum group) with $C(\IT
>\!\!\!\triangleleft \IZ_2) =C(ISO(\IT))$. \ethm
We can easily extend this result to higher dimensional commutative
tori, and can prove that the quantum isometry group coincides with
the classical isometry group. This is some kind of rigidity
result, and it will be interesting to investigate the nature of
quantum isometry groups of more general classical manifolds.

\noindent{\bf Example 2 : Noncommutative torus; holomorphic isomrtries}\\
Next we consider the simplest and well-known example of
noncommutative manifold, namely the noncommutative two-torus
$\cla_\theta$, where $\theta$ is a fixed irrational number (see
\cite{con}). It is the universal $C^*$ algebra generated by two
unitaries $U$ and $V$ satisfying the commutation relation
$UV=\lambda VU$, where $\lambda=e^{2 \pi i \theta}$. There is a
canonical faithful trace $\tau$ on $\cla_\theta$ given by
$\tau(U^m V^n)=\delta_{mn}$. We consider the canonical spectral
triple $({\cla^\infty}, \clh, D)$, where ${\cla^\infty}$ is the unital
 $\ast$-algebra spanned by
$U,V$,  $\clh=L^2(\tau) \oplus L^2(\tau)$ and $D$ is given by $$ D=\left( \begin{array}{cc} 0 &  d_1+id_2 \\
 d_1-id_2 &  0 \\ \end{array}\right),$$ where $d_1$ and $d_2$ are closed unbounded linear maps on $L^2(\tau)$
 given by $d_1(U^mV^n)=m U^mV^n,$ $d_2(U^mV^n)=nU^m V^n$. It is easy to compute the space of one-forms $\Omega^1_D$ (see \cite{jot}, \cite{fro}, \cite{con})
  and the Laplacian $\cll=- d^*d$ is given by $\cll(U^mV^n)=-(m^2+n^2)U^mV^n.$
  For simplicity of computation, instead of the full quantum isometry group we at first concentrate on an interesting quantum subgroup
  $\clg=QISO^{\rm hol}({\cla^\infty},\clh,D)$, which is the universal quantum group which leaves invariant the subalgebra of ${\cla^\infty}$
  consisting of polynomials in $U$, $V$ and $1$, i.e. span of $U^mV^n$ with $m,n \geq 0$.
  The proof of existence and uniqueness of such a universal quantum group is more or less identical to the proof of existence and uniqueness of QISO.
  We call $\clg$ the quantum group of ``holomorphic" isometries, and observe in the theorem
  stated below without proof (see \cite{jyotish}) that this
  quantum group is nothing but the quantum double torus studied in
  \cite{hajac}.
 \bthm
 Consider the following co-product $\Delta_\clb$ on the $C^*$ algebra
  $\clb=C(\IT^2) \oplus \cla_{2\theta}$, given on the generators $A_0,B_0,C_0,D_0$ as follows ( where $A_0,D_0$ correspond to $C(\IT^2)$ and $B_0,C_0$ correspond to $\cla_{2 \theta}$)
  $$ \Delta_\clb(A_0)=A_0 \ot A_0+C_0 \ot B_0, ~~\Delta_\clb(B_0)=B_0 \ot A_0+D_0 \ot B_0,$$
  $$ \Delta_\clb(C_0)=A_0 \ot C_0+C_0 \ot D_0,~~\Delta_\clb(D_0)=B_0 \ot C_0+D_0 \ot D_0.$$
  Then $(\clb,\Delta_0)$ is a compact quantum group and it has an action $\alpha_0$ on $\cla_\theta$ given by
  $$ \alpha_0(U)=U \ot A_0+V \ot B_0,~~\alpha_0(V)=U \ot C_0 + V \ot D_0.$$ Moreover, $(\clb,\Delta_\clb)$ is isomorphic (as quantum group) with
  $\clg=QISO^{\rm hol}({\cla^\infty},\clh, D)$.
 \ethm
 We refer to \cite{jyotish} for a proof of the above result, and to \cite{hajac} for
 the computation of the Haar stat and representation theory of the compact
  quantum group $\clg$.

\noindent{\bf Example 3 : Noncommutative Torus; full quantum
isometry group} By similar but somewhat tedious calculations (see
\cite{jyotish}) one can also describe explicitly the full quantum
isometry group $QISO({\cla^\infty},\clh,D)$. It is as a $C^*$ algebra has
eight direct summands, four of which are isomorphic with the
commutative algebra $C(\IT^2)$, and the other four are irrational
rotation algebras. \bthm $QISO(\cla_\theta)=\oplus_{k=1}^8
C^*(U_{k1},U_{k2})$ (as a $C^*$ algebra), where for odd $k$,
$U_{k1},U_{k2}$ are the two commuting unitary generators of
$C(\IT^2)$, and  for even $k$, $U_{k1}U_{k2}={\rm exp}(4 \pi i
\theta)U_{k2}U_{k1}$, i.e. they generate $\cla_{2 \theta}$. The
(co)-action on the generators $U,V$ (say) of $\cla_\theta$ are
given by the following : $$ \alpha_0(U)=U \ot (U_{11}+U_{31})+V
\ot (U_{52}+U_{62})+U^{-1} \ot (U_{21}+U_{41})+V^{-1} \ot
(U_{72}+U_{82}),$$ $$ \alpha_0(V)=U \ot (U_{51}+U_{71})+V \ot
(U_{12}+U_{22})+U^{-1} \ot (U_{61}+U_{81})+V^{-1} \ot
(U_{32}+U_{42}).$$ \ethm From the co-associativity condition, the
co-product of $QISO(\cla_\theta)$ can easily be calculated. For
the detailed description of the coproduct, counit, antipode and
study of the representation theory of $QISO(\cla_\theta)$, the
reader is referred to \cite{jyotish}. It is interesting to mention
here that the quantum isometry group of $\cla_\theta$ is a Rieffel
type deformation of the isometry group (which is same as the
quantum isometry group) of the commutative two-torus. The
commutative two-torus is a subgroup of its isometry group, but
when the isometry group is deformed into $QISO(\cla_\theta)$, the
subgroup relation is not respected, and the deformation of the
commutative torus, which is $\cla_{2\theta}$, sits in
$QISO(\cla_\theta)$ just as a $C^*$ subalgebra (in fact a direct
summand) but not as a quantum subgroup any more. This perhaps
provides some explanation of the non-existence of any Hopf algebra
structure on the
noncommutative torus. \vspace{2mm}\\
{\bf Acknowledgement :} The author would like to thank P. Hajac
for drawing his attention to the article \cite{hajac}, and S.L. Woronowicz for many valuable comments and suggestions which led to substantial improvement of the paper.

\end{document}